\documentclass[12pt]{article}

\usepackage{texdraw}
\usepackage{amsfonts}
\usepackage{amssymb}
\usepackage{latexsym}
\usepackage{epsf,graphicx}

\newcommand{\pl}{\partial}
\newcommand{\pr}{\noindent{\bf Proof.}\quad }
\newcommand{\epr}{\ $\blacksquare$ \vspace{2mm} }

\newcommand{\be} {\begin{eqnarray}}
\newcommand{\ee} {\end{eqnarray}}
\newcommand{\bep} {\begin{eqnarray*}}
\newcommand{\eep} {\end{eqnarray*}}

\newcommand {\Hol}{\mathop{\rm Hol}\nolimits}
\renewcommand {\Im}{\mathop{\rm Im}\nolimits}
\renewcommand {\Re}{\mathop{\rm Re}\nolimits}

\newcommand{\R}{{\mathbb R}}

\newcommand{\C}{{\mathbb C}}

\newtheorem{remar}{Remark}[section]
\newtheorem{examp}{Example}[section]

\newtheorem{corol}{Corollary}[section]

\newtheorem{theorem}{Theorem}[section]

\newcommand{\rema}{\begin{remar}\rm}
\newcommand{\erema}{\end{remar}}

\newcommand{\exa}{\begin{examp}\rm}
\newcommand{\eexa}{\end{examp}}

\def\lwvec(#1 #2){\linewd 0.1
           \lvec(#1 #2)
           \linewd 0.05}

\title{Boundary behavior and rigidity of semigroups of holomorphic mappings}

\author{Mark Elin
\\ {\small Department of Mathematics}
\\ {\small ORT  Braude College}
\\ {\small P.O. Box 78, 21982 Karmiel, Israel}
\\ {\small e-mail: mark$\_$elin@braude.ac.il}
\\ David Shoikhet
\\ {\small Department of Mathematics}
\\ {\small ORT Braude College}
\\ {\small P.O. Box 78, 21982 Karmiel, Israel}
\\ {\small e-mail: davs@braude.ac.il}}

\begin{document}

\maketitle

\begin{abstract}
In this paper we give some quantative characteristics of boundary
asymptotic behavior of semigroups of holomorphic self-mappings of
the unit disk including the limit curvature of their trajectories
at the boundary Denjoy--Wolff point. This enable us to establish
an asymptotic rigidity property for semigroups of parabolic type.

\vspace{4mm}

{\footnotesize Key words and phrases: holomorphic mapping,
asymptotic behavior, one-parameter semigroup, limit curvature,
rigidity.

}
\end{abstract}

\section{Introduction and main results}

We denote by $\Hol(D,\C)$ the set of all holomorphic functions on
a domain $D\subset\mathbb{C}$, and by $\Hol(D)$ the set of all
holomorphic self-mappings of $D$.

We say that a family $S=\left\{F_t\right\}_{t\geq
0}\subset\Hol(D)$ is \textit{a one-parameter continuous semigroup
on $D$} (semigroup, in short) if

(i) $F_{t}(F_{s}(z))=F_{t+s}(z)$ for all $t,s\geq 0$ and $z\in D$,

\noindent and

(ii) $\lim\limits_{t\rightarrow 0^+}F_{t}(z)=z$ for all $z\in D$.

In the case when $D$ is the open unit disk $\Delta=\{ z:|z|<1\}$,
it follows from a result of E. Berkson and H. Porta \cite{B-P}
that each semigroup is differentiable with respect to
$t\in\mathbb{R}^{+}=[0,\infty )$. So, for each one-parameter
continuous semigroup $S=\left\{F_{t}\right\}_{t\geq
0}\subset\Hol(\Delta)$, the limit
\[
\lim_{t\rightarrow 0^{+}}\frac{z-F_{t}(z)}{t}=f(z),\quad z\in
\Delta,
\]
exists and defines a holomorphic mapping
$f\in\Hol(\Delta,\mathbb{C})$. This mapping $f$ is called the
\textit{(infinitesimal) generator of} $S=
\left\{F_{t}\right\}_{t\geq 0}.$ Moreover, the function $
u(t,z):=F_{t}(z), \, (t,z)\in\mathbb{R}^{+}\times \Delta$, is the
unique solution of the Cauchy problem
\begin{equation}\label{2b}
\left\{
\begin{array}{l}
{\displaystyle\frac{\partial u(t,z)}{\partial
t}}+f(u(t,z))=0,\vspace{3mm}\\ u(0,z)=z,\quad z\in \Delta.
\end{array}
\right.
\end{equation}

In the same paper, Berkson and Porta proved that a function
$f\in\Hol(\Delta,\C)$ is a semigroup generator if and only if
there are a point $\tau\in\overline{\Delta}$ and a function
$p\in\Hol(\Delta,\C)$ with $\Re p(z)\ge0$, such that
\begin{equation}\label{bp}
f(z)=(z-\tau)(1-z\bar{\tau })p(z).
\end{equation}
This representation is unique. Moreover, if $S$ contains neither
the identity mapping nor an elliptic automorphism of $\Delta$,
then $\tau$ is a unique attractive fixed point of $S$, i.e.,
$\lim\limits_{t\to\infty}F_t(Z)=\tau,\ z\in\Delta,$ and
$\lim\limits_{r\to1^-}F_t(r\tau)=\tau$. The point $\tau$ is called
the \textit{Denjoy--Wolff point of} $S$.

In this paper we are interested in the boundary case
($\tau\in\pl\Delta$) for which the asymptotic behavior of the
semigroup has entirely different features then in the interior
case. It was shown in \cite{E-S1} that if $\tau\in\pl\Delta$, then
the angular derivative $\displaystyle
f'(\tau)=\angle\lim\limits_{z\to\tau}f'(z)=
\angle\lim\limits_{z\to1}\frac{f(z)}{z-\tau}$ of $f$ at the point
$\tau\in\pl\Delta$ exists and is a non-negative real number. One
distinguishes two cases: (a) $f'(\tau)>0$ (the \textit{hyperbolic
case}), and (b) $f'(\tau)=0$ (the \textit{ parabolic case}).

Although the asymptotic behavior of semigroups has been studied by
many mathematicians, the local geometry of semigroup trajectories
near the boundary Denjoy--Wolff point have been attracting an
intensive attention only recently.

It was shown that in many situations (in particular, always in the
hyperbolic case) the semigroup trajectories have tangent lines
passing through the Denjoy--Wolff point. Also, there is an
essential difference between hyperbolic and parabolic type
semigroups. Specifically, in the hyperbolic case the limit tangent
line depends on the initial point of the trajectory, while in the
parabolic case all the trajectories have the same tangent line (if
it exists). See Fig.~1 and \cite{C-DM, E-R-S-Y, E-S-Y,
E-K-R-S,E-S-book} for details.

\vspace{-25mm}

\begin{center}
\begin{figure}[h]
\centering
\includegraphics[angle=0,width=11cm,totalheight=8cm]{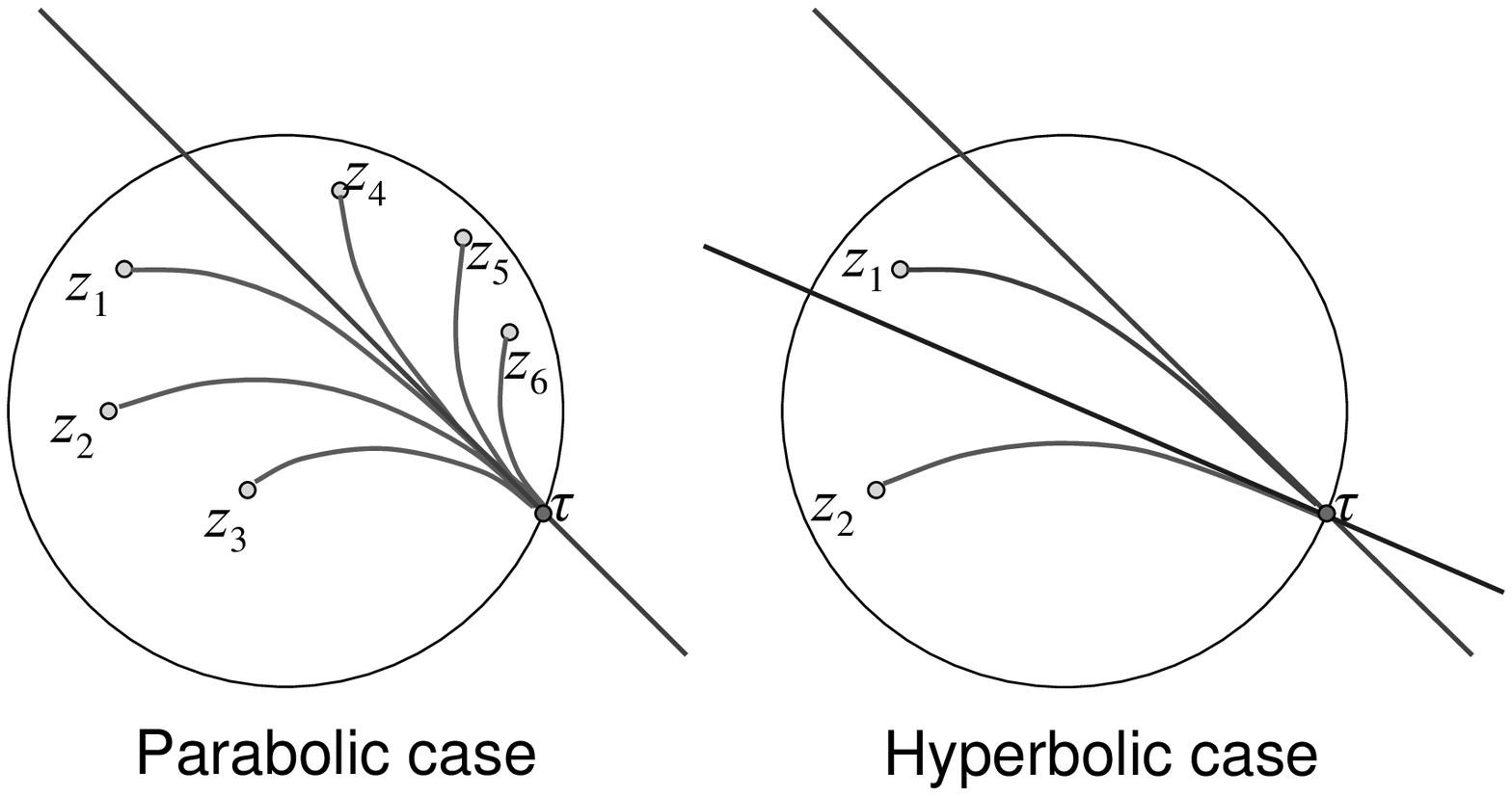}
\end{figure}
\end{center}

\vspace{-2.5cm}

An advanced question in this study is: how close is a semigroup
trajectory to its tangent line? In particular, one can ask: {\it
Is there a circle having the same tangent line, such that the
trajectories lie between this circle and the line?} See Fig.~2.

\begin{center}
\begin{figure}[h]
\centering
\includegraphics[angle=0,width=3.5cm,totalheight=4cm]{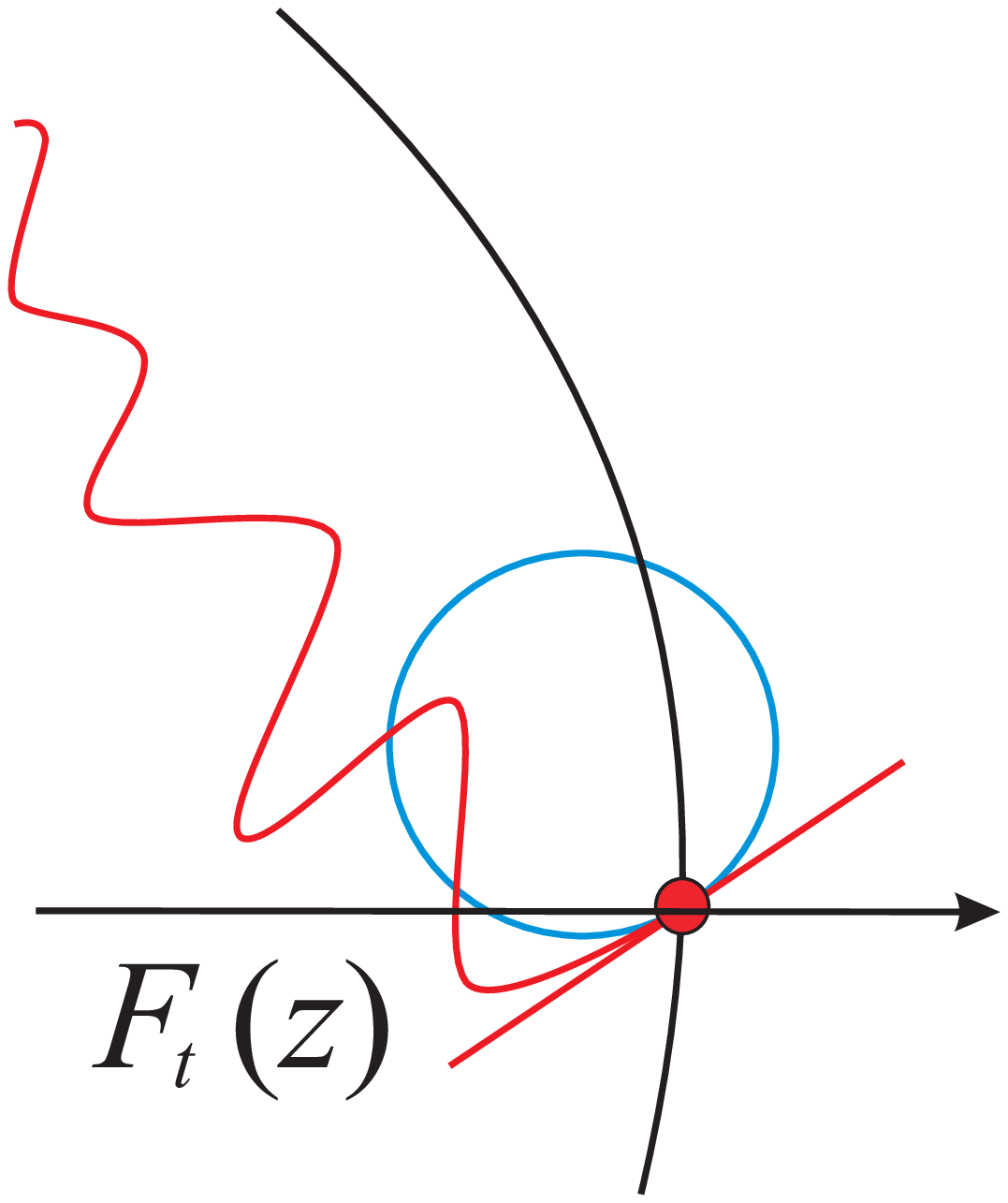}
\end{figure}
\end{center}

\vspace{-7mm}

Since any trajectory $\gamma_z=\left\{F_t(z),\ t\ge0 \right\},\
z\in\Delta,$ is an analytic curve, it has a finite curvature at
each its point $F_s(z)$. More precisely, given $z\in\Delta$ we
denote by $\kappa(z,s)$ the curvature of the trajectory $\gamma_z$
at the point $F_s(z)$ and by $\kappa(z)$ the \textit{limit
curvature of the trajectory}:
$\kappa(z):=\lim\limits_{s\to\infty}\kappa(z,s)$ (if it exists).
Therefore, the above question is equivalent to the following one:
{\it When the limit curvature of a semigroup trajectory exists
finitely?}

This question is closely connected with the so-called asymptotic
rigidity problem. Namely, one says that two semigroups with the
same Denjoy--Wolff point $\tau\in\pl\Delta$ have a similar
asymptotic behavior if
\[
\lim_{t\rightarrow \infty }\left( \arg \frac{1 -\bar\tau
F_{t}(z)}{1 -\bar\tau G_{t}(z)} \right) =0.
\]
From now on, we assume without loss of generality that $\tau=1$.

It turns out (see \cite{E-R-S-Y}) that two hyperbolic type
semigroups having similar asymptotic behavior actually coincide up
to rescaling. This fact is not longer true for parabolic type
semigroups. Moreover, if the generator $f$ of a semigroup
$S=\{F_t\}_{t\ge0}$ admits the representation
\[
f(z)=b(z-1)^2+o((z-1)^2),
\]
then for each $z\in\Delta$, the limit tangent line to the
trajectory $\gamma_z=\{{F_t(z),\ t\ge0}\}$ exists with
\[
\lim_{t\to\infty}\arg(1-F_t(z))=-\arg b,
\]
hence, does not depend on $z\in\Delta$ as well as on the remainder
$o((z-1)^2)$ (see \cite{E-R-S-Y} and Theorem~\ref{th-cur-par2}~(i)
below).

Therefore, a natural conjecture is: {\it If two semigroups having
similar asymptotic behavior and the same limit curvatures, then
they coincide up to rescaling.} The proof of this conjecture
follows directly from a more general fact (see
Theorem~\ref{th-rigid} below).

To answer the questions above for a semigroup $S=\{F_t\}_{t\ge0}$
generated by $f\in\Hol(\Delta,\C)$, we apply a linearization model
given by Abel's functional equation
\begin{equation}\label{abel}
h\left(F_t(z)\right)=h(z)+t.
\end{equation}
It is rather simple to see that the function $h:\Delta\mapsto\C$
defined by
\begin{equation}\label{h}
h'(z)=-\frac1{f(z)}\,,\quad h(0)=0,
\end{equation}
solves equation (\ref{abel}). This function is univalent and
convex in the positive direction of the real axis due to
(\ref{abel}). Sometimes it is called the {\it K{\oe}nigs function}
for the semigroup (see \cite{C-DM, E-R-S-Y, E-S-Y, Sis} and
\cite{E-S-book}).

It is more or less known that in the hyperbolic case each
semigroup trajectory $\gamma_z=\left\{F_t(z),\ t\ge0 \right\}$
converges to the boundary fixed point $\tau=1$ non-tangentially,
as $t\to\infty$. The following quantitative result in this
direction is given in \cite{E-R-S-Y}.

\begin{theorem}\label{teor4a}
Let $S=\{F_{t}\}_{t\geq 0}$ be a semigroup of hyperbolic type with
the Denjoy--Wolff point $\tau =1$ and let $f$ be its generator
with $f(1)=0$ and $a=f'(1)>0$. Then
\begin{equation} \label{A1}
\lim_{t\rightarrow \infty } \arg (1-F_{t}(z)) =a
\left[\lim_{r\rightarrow 1^-} \Im h(r) -\Im h(z)  \right].
\end{equation}
\end{theorem}

In the present paper we first complete Theorem~\ref{teor4a} by the
following assertion.

\begin{theorem}\label{th-cur-hyp1}
Let $S=\{F_{t}\}_{t\geq 0}$ be a semigroup of hyperbolic type with
the Denjoy--Wolff point $\tau =1$ and let $f$ be its generator.
Suppose that
\begin{equation}\label{c2-cond}
f(z)=a(z-1)+b(z-1)^2+R(z), \ \ \mbox{where }\
\angle\lim_{z\to1}\frac{R(z)}{(z-1)^2}=0.
\end{equation}
Then for any initial point $z_0\in\Delta$ the limit curvature
$\kappa(z_0)$ exists finitely. Moreover, the limit curvature
circle can be defined by the following equality
\begin{equation}\label{eq-hyp-cur1}
|1-z|^2\Im\bigl(b\overline{B}\bigr) +a\Im \bigl( (1-z)B \bigr)=0,
\end{equation}
where $B=\exp\left(a(h(z_0) - i \lim\limits_{r\to 1^-}\Im
h(r))\right)$. In particular, the limit curvature is zero if and
only if $\ \Im\bigl(b\overline{B}\bigr)=0$.
\end{theorem}

The parabolic case is more delicate. In this situation there are
semigroups which converge to the boundary Denjoy--Wolff point
tangentially as well as examples of non-tangentially converging
semigroups (see \cite{C-DM, C-M-P, E-K-R-S, E-R-S-Y}).

M. D. Contreras and S. D\'\i az-Madrigal in \cite{C-DM} have
considered the set $\mathrm{Slope}^+(\gamma_z)$ of all
accumulation points (as $t\to\infty$) of the function
$\arg\left(1-F_t(z)\right)$ and proved that these sets coincide
for all $z\in\Delta$. In addition, they have proven that if for a
function $h$ defined by (\ref{h}), the image $h(\Delta)$ lies in a
horizontal half-plane, then all the trajectories $\gamma_z$ tends
tangentially to $\tau=1$. In addition,
$\mathrm{Slope}^+(\gamma_z)$ is a single point which is equal to
either $\pi/Î2$ or $-\pi/2$. In general the question whether
$\mathrm{Slope}^+(\gamma_z)$ is a singlton is still open (see
\cite{C-DM, E-S-Y}).

Inasmuch as we are updated, all results known in this direction
require some smoothness conditions at the Denjoy--Wolff point. For
example, if the semigroup generator $f$ is twice differentiable at
the boundary Denjoy--Wolff point $\tau =1$, then all the
trajectories $\{F_{t}(z),\ t\ge0\}$ converge to this point
tangentially if and only if $\Re f''(1)=0$ (see \cite{E-R-S-Y}).

Furthermore, it may happen that for each $z\in \Delta$, there is a
horodisk $D(\tau ,k):=\left\{ \zeta \in \Delta :\ d(\zeta ,\tau
)<k\right\} ,\ k=k(z),$ internally tangent to the unit circle at
the point $\tau $, such that the trajectory $\{F_{t}(z)\}_{t\geq
0}$ lies outside $D(\tau ,k)$. In this case we say that the
semigroup $S$ converges to $\tau$ \textit{strongly tangentially}.
It is clear that the supremum of radii of such horodisks coincides
with the limit curvature radius. Conversely, if a semigroup
converges tangentially but not strongly tangentially, then its
trajectories have infinite limit curvature.

For three times differentiable generators at the boundary
Denjoy--Wolff point $\tau$ with $f'''(\tau)=0$ the following
rigidity phenomenon was established in \cite{SD-06}: {\it The
semigroup $S$ generated by $f$ converges to $\tau$ strongly
tangentially if and only if it consists of parabolic automorphisms
of $\Delta$ (i.e., its trajectories have finite curvature). }

Actually, the question on the finiteness of the limit curvature is
a general problem which is also relevant for non-tangentially
converging semigroups. In addition, one can ask: whether
$\kappa(z)$ might be finite for some points $z\in\Delta$ and
infinite to others? The following theorem answers these questions
for semigroups generated by functions which are
$(3+\varepsilon)$-smooth at the Denjoy--Wolff point.

\begin{theorem}\label{th-cur-par1}
Let $S=\left\{F_t \right\}_{t\ge0}$ be a semigroup generated by
$f\in\mathcal{C}^{3+\varepsilon}(1)$, i.e.,
\begin{equation}\label{gen1} f(z)=b(z-1)^2+c(z-1)^3+R(z),
\end{equation}
where $R\in\Hol(\Delta,\mathbb{C}),\ \displaystyle
\lim\limits_{z\to1}\frac{R(z)}{(z-1)^{3+\varepsilon}}=0$, and let
$b\not=0$.

(a) If $\displaystyle\ \Im\frac{c}{b^2}\not=0$, then all of the
trajectories 
have infinite limit curvature, i.e., $\kappa(z_0)=\infty$ for each
$z_0\in\Delta$.

(b) Otherwise, if $\displaystyle\ \Im\frac{c}{b^2}=0$, then each
trajectory $\{F_t(z_0),\ t\ge0\}$ has a finite limit curvature,
namely, $\kappa(z_0)=\left|\frac{2C}{b}\right|$, where
\[
C=|b|^2\Im h(z_0) + \Im b+ \int\limits_0^\infty
\Im\left(\frac{f\left(F_s(1)\right)\overline{b}}{(1-F_s(0))^2} -
\frac{c\overline{b}}{b(s+1)}\right)ds .
\]
\end{theorem}

Thus, under above assumptions if $\kappa(z)$ is finite for some
$z\in\Delta$, then it must be finite for all $z\in\Delta$.

Comparing Theorems~\ref{th-cur-hyp1} and \ref{th-cur-par1}, we see
again a cardinal difference between semigroups of hyperbolic and
parabolic types: in the hyperbolic case with some smoothness
conditions the limit curvature is always finite while in the
parabolic case the limit curvature may be infinite. At the same
time, it follows from \cite{SD-06} that if the second derivative
$f''(1)$ is purely imaginary, then the third derivative $f'''(1)$
should be real. So, an immediate consequence of part (b) of
Theorem~\ref{th-cur-par1} is the following fact.

\begin{corol}
Let $\left\{F_t\right\}_{t\ge0}$ be a semigroup of holomorphic
self-mappings of the open unit disk $\Delta$ generated by
$f\in\mathcal{C}^{3+\varepsilon}(1)$ of the form~(\ref{gen1}) with
$b\not=0$. If $\Re b=0$, then each semigroup trajectory converges
to $\tau=1$ strongly tangentially.
\end{corol}

As a matter of fact, Theorem~\ref{th-cur-par1} is based on the
following general result which contains complete quantitative
characteristics of the asymptotic behavior for semigroups
generated by functions smooth enough at the boundary Denjoy--Wolff
points.

\begin{theorem}\label{th-cur-par2}
Let $S=\left\{F_t\right\}_{t\ge0}$ be a continuous semigroup of
holomorphic self-mappings of the open unit disk $\Delta$ and let
$f$ be its generator.

(i) Suppose that $f$ admits the following representation:
\begin{equation}\label{circ_1}
f(z)=b(z-1)^2+R(z),
\end{equation}
where $R\in\Hol(\Delta,\mathbb{C}),\ \displaystyle
\lim\limits_{z\to1}\frac{R(z)}{(z-1)^2}=0$. Then
\begin{equation}\label{circ_2}
\frac1{1-F_t(z)}=-bt+G(z,t),\quad\mbox{where
}\lim_{t\to\infty}\frac{G(z,t)}t=0,
\end{equation}
and
\begin{equation}\label{circ_3}
\lim_{t\to\infty}
\left(\frac1{1-F_t(z)}-\frac1{1-F_t(0)}\right)=-bh(z).
\end{equation}

(ii) If $b\not=0$ and $R$ in (\ref{circ_1}) is of the form
$R(z)=c(z-1)^3+R_1(z)$ with $R_1\in\Hol(\Delta,\mathbb{C}),\
\displaystyle \lim\limits_{z\to1}\frac{R_1(z)}{(z-1)^3}=0$, i.e.,
$f$ admits the representation:
\begin{equation}\label{circ_4}
f(z)=b(z-1)^2+c(z-1)^3+R_1(z),
\end{equation}
then
\begin{eqnarray}\label{circ_5}
\frac1{1-F_t(z)}=-bt-\frac{c}{b}\log(t+1)+G_1(z,t),
\end{eqnarray}
where
$\displaystyle\lim_{t\to\infty}\frac{G_1(z,t)}{\log(t+1)}=0$, and
\begin{equation}\label{circ_6}
\lim_{t\to\infty} t\left(\frac1{1-F_t(z)}-\frac1{1-F_t(0)} + b
h(z) \right)= -\frac{c}{b}\, h(z).
\end{equation}

(iii) If function $R_1$ in (\ref{circ_4}) satisfies the condition
$\displaystyle\lim\limits_{z\to1}\frac{R_1(z)}{(z-1)^{3+\varepsilon}}=0$
for some $\varepsilon>0$, then there is a constant $A$ such that
\begin{equation}\label{circ_7}
\frac1{1-F_t(z)}=-bt-\frac{c}{b}\log(t+1) -b h(z)+A +G_2(z,t),
\end{equation}
where $\lim\limits_{t\to\infty}G_2(z,t)=0$.
\end{theorem}

Assertion (i) of Theorem~\ref{th-cur-par2} was proven in
\cite{E-R-S-Y}. For the sake of completeness we prove it below in
another way.

Now we are at the point to formulate our rigidity result which is
also a consequence of Theorem~\ref{th-cur-par2}.

\begin{theorem}\label{th-rigid}
Let $\left\{F_t \right\}_{t\ge0},\ \left\{G_t
\right\}_{t\ge0}\subset\Hol(\Delta)$ be two continuous semigroups
of holomorphic self-mappings of the open unit disk $\Delta$
generated by mapping $f$ and $g$, respectively. Suppose that $f$
and $g$ admit the following representations:
\begin{eqnarray*}
f(z)=b(z-1)^2+c_1(z-1)^3+r_1(z),\\
g(z)=b(z-1)^2+c_2(z-1)^3+r_2(z),
\end{eqnarray*}
where
$\displaystyle\lim\limits_{z\to1}\frac{r_1(z)}{(z-1)^{3+\varepsilon}}=
\lim\limits_{z\to1}\frac{r_2(z)}{(z-1)^{3+\varepsilon}}=0$. If
\[
\lim_{t\to\infty}\Im\left[\bar b\left(\frac1{1-F_t(z)} -
\frac1{1-G_t(z)} \right)\right]=0,
\]
then the semigroups coincide.
\end{theorem}

In the particular case when the limit curvature is finite, we
affirm the asymptotic rigidity conjecture mentioned above.

\begin{corol}
If the trajectories of two parabolic type semigroups generated by
mappings of the class $\mathcal{C}^{3+\varepsilon}(1)$, have the
same limit curvature circles, then the semigroups coincide up to
rescaling.
\end{corol}

\section{The right half-plane model}

\setcounter{equation}{0}

For some technical reasons we first transfer the study of the
semigroup behavior from the open unit disk to right half-plane by
using the Cayley transform
\[
C(z)=\frac{1+z}{1-z}\,.
\]
Now, given a semigroup $S=\left\{F_t\right\}_{t\geq
0}\subset\Hol(\Delta)$ with the Denjoy--Wolff point $\tau=1$, we
construct the semigroup $\Sigma=\left\{\Phi_t\right\}_{t\geq 0}$
of holomorphic self-mappings of the right half-plane
$\Pi=\left\{w\in\C:\ \Re w>0\right\}$ having the Denjoy--Wolff
point at $\infty$ as follows:
\begin{equation}\label{Phi-a}
\Phi_t(w)=C\circ F_t\circ C^{-1}(w).
\end{equation}
If $S$ is continuous (hence, differentiable) in $t$, then $\Sigma$
is too. More precisely, let $f$ be the infinitesimal generator of
$S$. Then by (\ref{bp}), $f$ must be of the form
$f(z)=-(1-z)^2p(z)$ with $\Re p(z)\ge0,\ z\in\Delta$.
Differentiating $\Phi_t$ given by (\ref{Phi-a}) at $t=0^+$, we
find that $\Sigma$ is generated by the mapping $-\phi$, where
\begin{equation}\label{phi}
\phi(w)=2p\left(C^{-1}(w)\right)
\end{equation}
(cf., \cite[Lemma 3.7.1]{E-R-S-2004}). So,
$\phi\in\Hol(\Pi,\overline{\Pi})$ and the semigroup
$\Sigma=\left\{\Phi_t\right\}_{t\geq 0}$ satisfies the Cauchy
problem
\begin{equation}\label{Cauch1}
\left\{
\begin{array}{l}
{\displaystyle\frac{\partial \Phi_t(w)}{\partial
t}}=\phi\left(\Phi_t(w)\right),\vspace{3mm}\\
\left.\Phi_t(w)\right|_{t=0}=w,\quad w\in \Pi.
\end{array}
\right.
\end{equation}

Concerning the K{\oe}nigs function $h$ defined by (\ref{h}), one
can modify it to $\sigma:=h\circ C^{-1}$. By direct calculations
we check that this modified function satisfies Abel's functional
equation
\begin{equation}\label{abel1}
\sigma\left(\Phi_t(w)\right)=\sigma(w)+t,\quad w\in\Pi,
\end{equation}
as well as the initial value problem:
\begin{equation}\label{h-1}
\sigma'(w)=\frac1{\phi(w)}\,,\quad \sigma(1)=0.
\end{equation}

It was already mentioned that the angular derivative $a=f'(1)$
always exists. Furthermore, by the Berkson--Porta
representation~(\ref{bp}) with $\tau=1$ and formula~(\ref{phi}),
the function $\phi$ can be represented as follows:
\begin{equation}\label{phi1}
\phi(w)=a(w+1)+\varrho(w),\quad\mbox{where}\quad
\angle\lim_{w\to\infty}\frac{\varrho(w)}{w}=0.
\end{equation}
If, in addition, $f\in \mathcal{C}^2(1)$, that is, $f$ admits
representation~(\ref{c2-cond}) with $a\ge0$, then
\begin{equation}\label{phi2}
\phi(w)=a(w+1)-2b+\varrho(w),\ \ \mbox{where}\
\lim_{w\to\infty}\varrho(w)=0.
\end{equation}

Suppose now that the semigroup $S$ generated by $f$ is of
parabolic type. Then we have $a=0$ in (\ref{c2-cond}). So,
formula~(\ref{phi2}) becomes
\begin{equation}\label{phi2-p}
\phi(w)=-2b+\varrho(w),\ \ \mbox{where}\
\lim_{w\to\infty}\varrho(w)=0.
\end{equation}

In the case when $S$ is of parabolic type and $f\in
\mathcal{C}^{3+\varepsilon}(1)$ for some $\varepsilon\ge0$, we
obtain in the same manner that
\begin{equation}\label{phi3-p}
\phi(w)=-2b+\frac{4c}{w+1}+\varrho(w),\ \ \mbox{where}\
\lim_{w\to\infty}(w+1)^{1+\varepsilon}\varrho(w)=0.
\end{equation}

Since $\Sigma$ has the Denjoy--Wolff point at $\infty$, we have by
Julia's Lemma (see, for example, \cite{SJH-93, SD, E-S-book}) that
$\Re \Phi_t(w)$ is an increasing function in $t\ge0$. The
tangential convergence of the semigroup means that the function
$\displaystyle \frac{\Im \Phi_t(w)}{\Re\Phi_t(w)}$ is unbounded as
$t$ tends to infinity. Roughly speaking, the semigroup converges
tangentially when $|\Im \Phi_t(w)|$ grows faster than $\Re
\Phi_t(w)$. Moreover, the original semigroup $S=\{F_{t}\}_{t\geq
0}$ converges strongly tangentially if and only if the function
$\Re\Phi_t(w)$ is bounded for each $w$ with $\Re w>0$. For this
reason, strongly tangentially convergent semigroups were referred
to in \cite{C-M-P} as semigroups of \textit{finite shift}; and
weakly tangentially convergent semigroups as semigroups of
\textit{infinite shift}. Therefore, a semigroup $S$ converges
strongly tangentially if and only if each trajectory of the
semigroup $\Sigma$ defined by (\ref{Phi-a}) has a vertical
asymptote. More generally, {\it a semigroup trajectory in the open
unit disk has a finite limit curvature if and only if the
corresponding trajectory in the right half-plane has an asymptote
as $t\to\infty$.}

Next we will consider hyperbolic and parabolic type semigroups in
the right half-plane separately.

\section{The hyperbolic case}

\setcounter{equation}{0}

\begin{theorem}\label{th-cur-hyp2}
Let $\{\Phi_{t}\}_{t\geq 0}\in\Hol(\Pi)$ be a semigroup of
hyperbolic type with the Denjoy--Wolff point at $\infty$ generated
by mapping $-\phi$ with $\Re\phi(w)\ge0,\ w\in\Pi$. Suppose that
\begin{equation}\label{repr_c}
\phi(w)=\alpha w+\beta+\varrho(w), \ \ \mbox{where }\
\angle\lim_{w\to\infty}\varrho(w)=0.
\end{equation}
Then for any initial point $w_0\in\Pi$, the trajectory
$\left\{\Phi_s(w_0),\ s\ge0\right\}$ has the asymptote defined by
the following equality
\begin{eqnarray}\label{line5-c}
\Im\left[\left(\alpha w+\beta\right)\overline{B}\right] = 0,
\end{eqnarray}
where $B=\exp\left(\alpha(\sigma(w_0) - i
\lim\limits_{x\to\infty}\Im\sigma(x))\right)$ and $\sigma$ is
defined by (\ref{h-1}).
\end{theorem}

\pr Consider the semigroup $S=\left\{F_t\right\}_{t\ge0}$ in the
open unit disk defined by $F_t=C^{-1}\circ\Phi_t\circ C$.
Comparing formulas~(\ref{c2-cond}), (\ref{phi2}) and
(\ref{repr_c}) we conclude that $S$ is generated by the mapping
$f$ of the form (\ref{c2-cond}) with $a=\alpha\quad\mbox{and}\quad
b=\frac{\alpha-\beta}2.$ It was shown in \cite[Theorem 1 and
Remark 3]{E-L-R-S} that in this case for each $t\ge0$, the
semigroup element $F_t\in\Hol(\Delta)$ is twice differentiable at
the point $\tau=1$ and
\[
F_t'(1)=e^{-t\alpha},\quad F_t''(1)=\frac{\alpha-\beta}\alpha
e^{-t\alpha}(e^{-t\alpha}-1).
\]
Using (\ref{Phi-a}), a direct calculation  shows that
\begin{equation}\label{Phi-b}
\Phi_t(w)=e^{t\alpha}w
+\frac\beta\alpha(e^{t\alpha}-1)+\Gamma(t,w),
\end{equation}
where $\lim\limits_{w\to\infty}\Gamma(t,w)=0$. In particular, for
each $w\in\Pi$ we have:
\begin{equation}\label{Phi-c}
\lim\limits_{w\to\infty}\frac{\Phi_t(w)}w=e^{t\alpha}.
\end{equation}

Furthermore, by (\ref{Phi-a}) we have
\[
\lim\limits_{t\to\infty}\arg\Phi_t(w)=-
\lim\limits_{t\to\infty}\arg\left(1-F_t(C^{-1}(w)) \right).
\]
Therefore, formula (\ref{A1}) in Theorem~\ref{teor4a} can be
rewritten as
\begin{equation}\label{N}
\lim\limits_{t\to\infty}\arg\Phi_t(w)=\alpha \left[\Im\sigma(w) -
\lim\limits_{x\to+\infty}\Im\sigma(x)\right].
\end{equation}
Denote
\begin{equation}\label{N+1}
g_t(w)=\frac{\Phi_t(w)}{\left|\Phi_t(1)\right|}\,.
\end{equation}
It follows by the continuous version of the Valiron theorem (see,
for example, \cite{Val, E-S-V}), that the limit
$g(w)=\lim\limits_{t\to\infty}g_t(w)$ exists for all $w\in\Pi$ and
$g$ satisfies Schr\"oder's functional equation
\[
g\left(\Phi_s(w)\right)=e^{s\alpha}g(w),\quad s\ge0,
\]
which is actually equivalent to the differential equation
\begin{equation}\label{N+2}
g'(w)\phi(w)=\alpha g(w)
\end{equation}
(cf., \cite{SD, Sis}). In addition,
\[
\lim\limits_{t\to\infty}\arg\Phi_t(w)=\arg g(w)
\]
and $g(1)=e^{-i\theta}$ for some $\theta\in\R$. Since
$\sigma(1)=0$ (see formula~(\ref{h-1})) and
$g(w)=\lim\limits_{t\to\infty}g_t(w)$, one calculates $\theta$ by
using~(\ref{N}) and~(\ref{N+1}). Namely,
\begin{eqnarray*}
\theta=-\arg g(1)=-\lim_{t\to\infty}\arg g_t(1)
=-\lim_{t\to\infty}\arg \Phi_t(1) \\ = -\alpha\left[
\Im\sigma(1)-\lim_{x\to\infty}\Im\sigma(x) \right]
=\alpha\lim_{x\to\infty}\Im\sigma(x).
\end{eqnarray*}

To proceed we denote by $L_t(w_0)$ the tangent line to the
trajectory $\left\{\Phi_s(w_0),\ s\ge0\right\}$ at the point
$\Phi_t(w_0)$ for a fixed $w_0\in\Pi$. Its equation is
\begin{equation}
\Im \left[\biggl(w-\Phi_t(w_0)\biggr)
\overline{\phi(\Phi_t(w_0))}\right]=0,
\end{equation}
or, what is one and the same,
\begin{equation}\label{as1}
\Im\left[ w\overline{\phi\left(\Phi_t(w_0)\right)}\right] =
\Im\left[ \Phi_t(w_0)\overline{\phi(\Phi_t(w_0))} \right].
\end{equation}
Substituting (\ref{repr_c}) into (\ref{as1}) we see that
\begin{equation}\label{line2-c}
\Im\left[w\overline{\left(\alpha\Phi_t(w_0)+\beta +
\varrho(\Phi_t(w_0))\right)}\right]=\Im\left[
\Phi_t(w_0)\overline{(\beta+\varrho(\Phi_t(w_0)))}\right].
\end{equation}

Dividing now equation~(\ref{line2-c}) by $\left|\Phi_t(1)\right|$,
we have
\begin{eqnarray}\label{line3-c}
\Im\left[w\cdot\overline{\left(\alpha
g_t(w_0)+\frac{\beta}{\left|\Phi_t(1)\right|}
+\frac{\varrho(\Phi_t(w_0))}{\left|\Phi_t(1)\right|}\right)}\right]
\nonumber \\= \Im\left[
g_t(w_0)\overline{(\beta+\varrho(\Phi_t(w_0)))}\right].
\end{eqnarray}
Letting $t$ tend to infinity, we conclude that the net of lines
$\left\{L_t(w_0)\right\}$ converges to the limit line $L(w_0)$ the
equation of which by~(\ref{line3-c}) is
\[
\alpha\Im\left[w\overline{g(w_0)}\right] = \Im\left[
g(w_0)\overline{\beta}\right],
\]
or
\[
\Im\left[\left(\alpha w+\beta \right)\overline{g(w_0)}\right] = 0,
\]
Comparing now (\ref{h-1}) with (\ref{N+2}), we get
$$g(w)=\exp\left(\alpha(\sigma(w)- i\lim\limits
_{x\to\infty}\Im\sigma(x))\right).$$ Setting $B=g(w_0)$, we
complete the proof. \epr

\noindent{\bf Proof of Theorem \ref{th-cur-hyp1}.} Suppose that
the generator $f$ of a semigroup $S$ is of the form
(\ref{c2-cond}). Construct the semigroup $\Sigma$ in the right
half-plane and denote its generator by $-\phi$. Then $\phi$ has
form (\ref{repr_c}) with $\alpha=a$ and $\beta=a-2b$. We have
already shown that each semigroup trajectory $\left\{\Phi_t(w_0),\
t\ge0 \right\}$ has the asymptote defined by (\ref{line5-c}).
Substituting in this formula $\alpha,\beta$ and
$w=\frac{1+z}{1-z}$, we attain at (\ref{eq-hyp-cur1}). \epr

\section{The parabolic case}

\setcounter{equation}{0}

\begin{theorem}\label{th-cur-par3}
Let $\{\Phi_{t}\}_{t\geq 0}\in\Hol(\Pi)$ be a semigroup of
parabolic type with the Denjoy--Wolff point at $\infty$ generated
by mapping $-\phi$.

(i) Suppose that $\phi$ admits the representation
\begin{equation}\label{repr_cp}
\phi(w)=\beta+\varrho(w),
\end{equation}
where $\varrho\in\Hol(\Pi,\mathbb{C}),\
\lim\limits_{z\to\infty}\varrho(w)=0$. Then
\begin{equation}\label{arg}
\Phi_t(w)=\beta t+\Gamma(w,t),\quad\mbox{where
}\lim_{t\to\infty}\frac{\Gamma(w,t)}t=0,
\end{equation}
and
\begin{equation}\label{difference1}
\lim_{t\to\infty} \left(\Phi_t(w)-\Phi_t(1)\right)=\beta\sigma(w).
\end{equation}

(ii) If $\beta\not=0$ and $\varrho$ in (\ref{repr_cp}) is of the
form $\varrho(w)=\frac\gamma{w+1}+\varrho_1(w)$ with
$\varrho_1\in\Hol(\Pi,\C),\
\lim\limits_{w\to\infty}w\varrho(w)=0$, i.e., $\phi$ admits the
representation
\begin{equation}\label{repr_cp_new}
\phi(w)=\beta+\frac{\gamma}{w+1}+\varrho_1(w),
\end{equation}
then
\begin{equation}\label{repr_cp_new1}
\Phi_t(w)=\beta t+\frac\gamma\beta\log(t+1)+\Gamma_1(w,t),
\end{equation}
where $\lim\limits_{t\to\infty}\frac{\Gamma_1(w,t)}{\log(t+1)}=0,$
and
\begin{equation}\label{difference2}
\lim_{t\to\infty} t\left(\Phi_t(w)-\Phi_t(1) -
\beta\sigma(w)\right)=\frac{\gamma\sigma(w)}\beta\,.
\end{equation}

(iii) If function $\varrho_1$ in (\ref{repr_cp_new}) satisfies the
condition
$\lim\limits_{w\to\infty}w^{1+\varepsilon}\varrho_1(w)=0$ for some
$\varepsilon>0$, then there is a constant $A$ such that for all
$w\in\Pi$
\begin{equation}\label{repr_cp_new2}
\Phi_t(w)=\beta
t+\frac\gamma\beta\log(t+1)+\beta\sigma(w)+A+\Gamma_2(w,t),
\end{equation}
where $\lim\limits_{t\to\infty}\Gamma_2(w,t)=0$.
\end{theorem}

\pr (i) Fix $w\in\Pi$ and consider $\Phi_t(w)$ as a (complex
valued) function of the real variable~$t$. This function tends to
infinity as $t\to\infty$. Thus, by the L'Hospital rule
\[
\lim_{t\to\infty}\frac{\Phi_t(w)}t=\lim_{t\to\infty}\frac{\phi\left(\Phi_t(w)\right)}1=\beta.
\]
Furthermore,
\begin{eqnarray*}
\lim_{t\to\infty} \left(\Phi_t(w)-\Phi_t(1)\right)=
\lim_{t\to\infty}\int_1^w\left(\Phi_t(z)\right)'dz =
\lim_{t\to\infty}\int_1^w \frac{\phi\left(\Phi_t(z)
\right)}{\phi(z)}dz
\\= \lim_{t\to\infty}\int_1^w  \frac{\beta+\varrho\left(\Phi_t(z)
\right)}{\phi(z)}dz =\beta\int_1^w  \frac{dz}{\phi(z)} =
\beta\sigma(w).
\end{eqnarray*}

(ii) To prove this assertion we need to show that
\begin{eqnarray}\label{limit1}
\lim_{t\to\infty}\frac1{\log(t+1)}\cdot\left(\Phi_t(w)-\beta t-
\frac\gamma\beta\log(t+1) \right)=0.
\end{eqnarray}
To do this, we calculate
\begin{eqnarray*}
&&\Phi_t(w)-\beta t-\frac\gamma\beta \log(t+1)
=\int_0^t\left(\Phi_s(w)-\beta s-\frac\gamma\beta\log(s+1)
\right)'ds +w \\ && = \int_0^t
\left(\phi\left(\Phi_s(w)\right)-\beta-\frac\gamma{\beta(s+1)}
\right)ds +w \\ && = \int_0^t
\left(\frac\gamma{\Phi_s(w)+1}-\frac\gamma{\beta(s+1)} +
\varrho_1\left(\Phi_s(w)\right) \right)ds +w.
\end{eqnarray*}
Note that
\begin{eqnarray*}
&&\frac\gamma{\Phi_s(w)+1}-\frac\gamma{\beta(s+1)} + \varrho_1
\left( \Phi_s(w)\right)\\ && = \frac1{s+1}\cdot\left(
\frac{\gamma(s+1)}{\Phi_s(w)+1}-\frac\gamma\beta+
\frac{\Phi_s(w)\varrho_1 \left(\Phi_s(w)\right)}{\Phi_s(w)/(s+1)}
\right),
\end{eqnarray*}
where by (i) and our assumption the second factor tends to zero.
Therefore, for each $\varepsilon>0$ there is $t_0$ such that
\[
\left| \frac{\gamma(s+1)}{\Phi_s(w)+1}-\frac\gamma\beta+
\frac{\Phi_s(w)\varrho\left(\Phi_s(w)\right)}{\Phi_s(w)/(s+1)}
\right|<\frac\varepsilon2
\]
for $s>t_0$. Hence
\begin{eqnarray*}
&& \left|\int_0^t
\left(\frac\gamma{\Phi_s(w)+1}-\frac\gamma{\beta(s+1)} +
\varrho_1\left(\Phi_s(w)\right) \right) ds\right|\le  \\ \le &&
\left|\int_0^{t_0}
\left(\frac\gamma{\Phi_s(w)+1}-\frac\gamma{\beta(s+1)} +
\varrho_1\left(\Phi_s(w)\right) \right) ds\right| +  \int_{t_0}^t
\frac{\varepsilon/2}{s+1}  ds  .
\end{eqnarray*}
Now for $t_0$ given above, choose $t_1$ such that
\[
\left|\int_0^{t_0}
\left(\frac\gamma{\Phi_s(w)+1}-\frac\gamma{\beta(s+1)} + \varrho_1
\left(\Phi_s(w)\right) \right) ds\right|\le
\frac\varepsilon2\log(t_1+1).
\]
So, for $t\ge\max(t_0,t_1)$ we have
\begin{eqnarray*}
&& \left|\int_0^t \left(\frac\gamma{\Phi_s(w)+1} -
\frac\gamma{\beta(s+1)} + \varrho_1\left(\Phi_s(w)\right)
\right)ds \right| \le \\ && \le \frac\varepsilon2\log(t_1+1) +
\int_{t_0}^t \frac{\varepsilon/2}{s+1}  ds  \le
\varepsilon\log(t+1).
\end{eqnarray*}

Since $\varepsilon>0$ can be chosen arbitrarily, we conclude that
the limit in (\ref{limit1}), hence, in (\ref{repr_cp_new1}), is
zero. In addition, we have that
\begin{eqnarray*}
&&\lim_{t\to\infty} t\left(\Phi_t(w)-\Phi_t(1)
-\beta\sigma(w)\right) = \lim_{t\to\infty} t\cdot \int_1^w
\left(\Phi_t(z)-\beta\sigma(z)\right)'dz
\\ = && \lim_{t\to\infty} t\cdot \int_1^w
\frac{\phi\left(\Phi_t(z)\right)-\beta}{\phi(z)} dz =
\lim_{t\to\infty} t \cdot \int_1^w
\frac{\frac\gamma{\Phi_t(z)}+\varrho_1\left(\Phi_t(z)\right)}{\phi(z)}dz
\\ = &&  \lim_{t\to\infty} \int_1^w \frac{t}{\Phi_t(z)}\cdot
\frac{\gamma+\Phi_t(z)\varrho_1\left(\Phi_t(z)\right)}{\phi(z)}dz
= \frac{\gamma\sigma(w)}\beta\,.
\end{eqnarray*}
This proves (\ref{difference2}).

(iii) First, we prove that for each $w\in\Pi$ the limit
\[
H(w):=\lim_{t\to\infty}\left(\Phi_t(w)-\beta
t-\frac\gamma\beta\log(t+1) \right)
\]
exists. Indeed,
\begin{eqnarray*}
&& \Phi_t(w)-\beta t-\frac\gamma\beta\log(t+1) \\ && = \int_0^t
\frac{d}{ds}\left(\Phi_s(w)- \beta s-\frac\gamma\beta\log(s+1)\right)ds +w \\
&& = \int_0^t \left(\phi\left(\Phi_s(w)\right)- \beta-
\frac\gamma{\beta(s+1)}\right)ds +w \\
&&  = \int_0^t \varrho_1\left(\Phi_s(w)\right)ds +
\frac\gamma\beta\int_0^t\left(\frac\beta{\Phi_s(w)+1}-\frac1{s+1}
\right)ds +w.
\end{eqnarray*}

To this end, we need to show that the two integrals
\[
\int_0^\infty \varrho_1\left(\Phi_s(w)\right)ds\quad
\mbox{and}\quad
\int_0^\infty\left(\frac\beta{\Phi_s(w)+1}-\frac1{s+1} \right)ds
\]
converge. Writing the first integral in the form
\[
\int_0^\infty \frac{\Phi_s^{1+\varepsilon}(w)\varrho_1
\left(\Phi_s(w)\right)}{\Phi_s^{1+\varepsilon}(w)}ds,
\]
we see that it converges since the nominator tends to zero and the
denominator behaves as $s^{1+\varepsilon}$ when $s\to\infty$ by
(\ref{arg}).

Consider the second integral. Using (\ref{repr_cp_new1}), we have:
\begin{eqnarray*}
&& \int_0^\infty \frac{\beta(s+1)-\Phi_s(w)-1}{(\Phi_s(w)+1)(s+1)} ds \\
&& = \int_0^\infty \frac{\beta-1 - \frac\gamma\beta\log(s+1)-
\Gamma_1(w,s)} {\left(\beta s+
\frac\gamma\beta\log(s+1)+\Gamma(w,s) +1 \right)(s+1)}\, ds  \\
&& = \int_0^\infty \frac{1}{(s+1)^2}\cdot \frac{\beta-1 -
\frac\gamma\beta\log(s+1)- \Gamma_1(w,s)}
{\beta-\frac{1-\beta}{s+1}+\frac{\gamma\log(s+1)}{\beta(s+1)}+\frac{\Gamma_1(w,s)}{s+1}}\,
ds \,.
\end{eqnarray*}
Since the second factor in the integrand is $O(\log(s+1))$, this
integral also converges.

Finally, by (\ref{difference1}) we get
\[
H(w)-H(1)=\lim_{t\to\infty}\left(\Phi_t(w)-\Phi_t(1)\right)=\beta\sigma(w).
\]
This proves formula (\ref{repr_cp_new2}) with $A=H(1)$ and
$\Gamma_2(w,t)=\Phi_t(w)-\beta t -\frac\gamma\beta\log(t+1)
-\beta\sigma(w)-A$. The proof of the theorem is complete. \epr

\begin{theorem}\label{th-cur-par4}
Let $\{\Phi_{t}\}_{t\geq 0}\in\Hol(\Pi)$ be a semigroup of
parabolic type generated by mapping $-\phi$, where $\phi$ has the
form
\[
\phi(w)=\beta+\frac\gamma{w+1}+\varrho_1(w)
\]
with $\varrho_1\in\Hol(\Pi,\mathbb{C}),\ \beta,\gamma\in\C$.
Assume that $\beta\not=0$ and $\lim\limits_{z\to\infty}
w^{1+\varepsilon}\varrho_1(w)=0$ for some $\varepsilon>0$. The
following assertions hold.

(a) If $\displaystyle\ \Im\frac\gamma{\beta^2}\not=0$, then for
any trajectory $\{\Phi_t(w),\ t\ge0\}$ there is no asymptote.

(b) Otherwise, if $\displaystyle\ \Im\frac\gamma{\beta^2}=0$, then
each trajectory $\{\Phi_t(w),\ t\ge0\}$ has an asymptote. In this
case, the asymptote equation is $\Im\left[w\cdot\bar
\beta\right]=B$ with
\[
B=|\beta|^2\Im \sigma(w) +  \int\limits_0^\infty
\Im\left(\phi\left(\Phi_s(1)\right)\bar{\beta} -
\frac{\gamma\bar{\beta}}{\beta(s+1)}\right)ds -\Im \beta .
\]
\end{theorem}

\pr It follows by Theorem~\ref{th-cur-par3}~(i) that
\[
\lim_{t\to\infty}\arg \Phi_t(w)= \lim_{t\to\infty}\arg
\frac{\Phi_t(w)}t=\arg \beta.
\]
So, if an asymptote exists, its equation must be
$\Im\left[w\cdot\bar \beta\right]=const$.

(a) If $\displaystyle\Im\frac\gamma{\beta^2}\not=0$, then by
(\ref{repr_cp_new1}) we have
\begin{eqnarray*}
\Im\left[\Phi_t(w)\cdot\bar \beta\right] = \Im\left[\left(\beta
t+\frac\gamma\beta\log(t+1)+\Gamma(w,t) \right)\cdot\bar \beta\right] \\
=\log(t+1)|\beta|^2\Im\frac\gamma{\beta^2}+ \Im\left[\Gamma(w,t)
\cdot\bar \beta\right].
\end{eqnarray*}
Thus, under our assumption we have that
$\Im\left[\Phi_t(w)\cdot\bar \beta\right]$ is unbounded, so there
is no asymptote for the trajectory $\left\{\Phi_t,\ t\ge0
\right\},\ w\in\Pi$.

(b) If $\displaystyle\Im\frac\gamma{\beta^2}=0$, then by
(\ref{repr_cp_new2}) we obtain that
\[
\Im\left[\Phi_t(w)\cdot\bar \beta\right] = \Im\left[\left(\beta
t+\frac\gamma\beta\log(t+1)+\beta\sigma(w)+A+\Gamma_2(w,t)\right)\cdot\bar
\beta\right],
\]
where $\lim\limits_{t\to\infty}\Gamma_2(w,t)=0$. Hence,
\[
\lim_{t\to\infty}\Im\left[\Phi_t(w)\cdot\bar \beta\right] =
|\beta|^2\Im \sigma(w) +\Im(A\bar\beta).
\]
So, in this case the asymptote for the trajectory
$\left\{\Phi_t(w),\ t\ge0\right\}$ does exist. But we have already
seen in the proof of Theorem~\ref{th-cur-par3} that
\[
A=H(1)= \int\limits_0^\infty \left(\phi\left(\Phi_s(1)\right)-
\beta - \frac\gamma{\beta(s+1)}\right)ds +1 ,
\]
and this completes our proof. \epr

Again we observe that Theorems~\ref{th-cur-par1}
and~\ref{th-cur-par2} are interpretations of
Theorems~\ref{th-cur-par4} and~\ref{th-cur-par3}, respectively,
for the semigroup $S=\{F_t \}_{t\ge0},\ F_t=C^{-1}\circ
\Phi_t\circ C,$ acting in the open unit disk.

\vspace{2mm}

\noindent{\bf Proof of Theorem~\ref{th-rigid}}. Assume that
\[
\lim_{t\to\infty}\Im\left[\bar b\left(\frac1{1-F_t(z)} -
\frac1{1-G_t(z)} \right)\right]=0.
\]
Using Theorem~\ref{th-cur-par2}~(iii), we have
\begin{eqnarray*}
&&\lim_{t\to\infty}\Im\left[\bar b\left(\frac1{1-F_t(z)} -
\frac1{1-G_t(z)} \right)\right]= \\ && \lim_{t\to\infty}\Im\left[
\frac{(c_2-c_1)\bar b}b\log(t+1) + |b|^2\left(h_2(z)-h_1(z)\right)
+ \bar b(A_1-A_2)\right],
\end{eqnarray*}
where $A_1, A_2$ are constants and $h_1,h_2$ are K{\oe}nigs
functions for $\left\{F_t \right\}_{t\ge0}$ and $\left\{G_t
\right\}_{t\ge0}$, respectively. The last limit may be zero only
if $\ {\Im\frac{(c_2-c_1)\bar b}b=0}$ and then $\Im
\left(h_1(z)-h_2(z)\right)=\Im\frac{A_1-A_2}b$. Therefore,
${h_1(z)-h_2(z)}=const.$ Since $h_1(0)=h_2(0)=0,$ we have $h_1(z)=
h_2(z)$. Now by (\ref{h}), we conclude that $f\equiv g$. \epr

For the sake of completeness we provide an equivalent result for
the right half-plane case.

\begin{theorem}\label{th-rigid1}
Let $\left\{\Phi_t \right\}_{t\ge0},\ \left\{\Psi_t
\right\}_{t\ge0}\subset\Hol(\Pi)$ be two continuous semigroups of
holomorphic self-mappings in the right half-plane generated by
mappings $-\phi$ and $-\psi$, respectively. Suppose that $\phi$
and $\psi$ admit the following representations:
\begin{eqnarray*}
\phi(w)=\beta+\frac\gamma{w+1}+\varrho(w),\\
\psi(w)=\beta+\frac{\gamma_1}{w+1}+\varrho_1(w),
\end{eqnarray*}
with $\varrho,\varrho_1\in\Hol(\Pi,\mathbb{C}),\
\lim\limits_{z\to\infty}w^{1+\varepsilon}\varrho(w)=
\lim\limits_{z\to\infty}w^{1+\varepsilon}\varrho_1(w)=0$,
$\beta,\gamma,\gamma_1\in\C$, and let $\beta\not=0$. If
\[
\lim_{t\to\infty}\Im\left[\bar \beta\left(\Phi_t(w) - \Psi_t(w)
\right)\right]=0,
\]
then the semigroups coincide.
\end{theorem}

{\bf Acknowledgments.} This research is part of the European
Science Foundation Networking Programme HCAA. The second author is
supported by the European Community project TODEQ
(MTKD-CT-2005-030042).


\end{document}